\documentclass{amsart}

\usepackage{amssymb,enumerate,graphics}
\usepackage{collection} 

\usepackage{array}
\usepackage{rotating}

\usepackage[plain]{algorithm} % possible options: plain, ruled, boxed
\usepackage{algorithmic}
\algsetup{indent=2em}
\floatname{algorithm}{\sc Algorithm} % make it AMSTeX like
\renewcommand{\algorithmicrequire}{\textbf{Input:}}
\renewcommand{\algorithmicensure}{\textbf{Output:}}
\renewcommand{\algorithmiccomment}[1]{\hfill #1}

\newcommand{\m}{\bar}

\newtheorem{Thm}{Theorem}[section]

\newtheorem{Prop}[Thm]{Proposition}

\title[Computing in unipotent groups]
{Computing in \\unipotent and reductive algebraic groups}

\date{\today}

\author{Arjeh M. Cohen \and Sergei Haller \and Scott H.\ Murray}

\thanks{This paper was written during a stay of the
second author at the Magma group at University of Sydney.
The authors would like to thank Bj\"orn Assmann, John Cannon, Michael Harrison, 
Anthony Henderson, Bill Kantor, Allan Steel, and Donald E.\ Taylor for useful 
suggestions and discussions.}

\address{
Department of Mathematics and Computer Science\\
Eindhoven University of Technology\\
PO Box 513\\
5600 MB Eindhoven\\
Netherlands}
\email{\eml{A.M.Cohen@tue.nl}}

\address{
School of Mathematics and Statistics F07 \\
Faculty of Science \\
University of Sydney, NSW 2006 \\
Australia}
\email{\eml{sergei@sergei-haller.de},
       \eml{murray@maths.usyd.edu.au}}

\begin{document}

\begin{abstract}
The unipotent groups are an important class of algebraic groups.
We show that techniques used to compute with finitely generated nilpotent groups 
carry over to unipotent groups.
We concentrate particularly on the maximal unipotent subgroup of a 
split reductive group and show how this improves computation in the
reductive group itself.
\end{abstract}

\maketitle

\section{Introduction}\label{S-intro}
A linear algebraic group is \emph{unipotent} if its elements are unipotent (ie, 
the elements only have eigenvalue one in every representation).
Unipotent groups play a prominent role in the theory of algebraic groups.
In this paper, we present algorithms for efficient element operations in
unipotent groups.
Since unipotent groups are nilpotent, we adapt methods used for computing
in finitely generated nilpotent groups.
A \emph{PC group} provides a unique 
computer representation for the elements of nilpotent groups (see, for example,
\cite{HoltEickOBrien05Book}).
\emph{Collection} gives efficient algorithms for multiplication and inversion of
group elements.
We modify these concepts to work with a large class of unipotent
groups defined over a field.
This class contains all unipotent groups if the field has characteristic zero.
It also contains the full unipotent subgroup of every split reductive group.

Algorithms for element operations in split reductive algebraic groups are given 
in \cite{CohenMurrayTaylor04}.
Computations in the unipotent subgroup make the largest single contribution to the time 
taken by these algorithms. 
This was the main impetus for the current paper, in which we prove:
\begin{Thm}\label{T-main}
Let $\F$ be a field with effective algorithms for the basic element operations.
Let $G$ be a split reductive algebraic group over $\F$ with rank $n$.
Then there is a normal form for elements of $G(\F)$.
The word problem for elements in normal form requires $O(n^2)$ field operations, 
and multiplying or inverting them requires $O(n^3)$ field operations. 
\end{Thm}
\noindent This theorem is a great improvement over the analysis of
\cite{CohenMurrayTaylor04}, where we proved that the operations are polynomial
time, but did not compute the exponent.
This result is optimal in the sense that the timings are asymptotically the same 
as the straightforward methods for matrices.

In Section~\ref{S-pres}, we construct \emph{FC group schemes}, which give normal
forms for elements of unipotent groups.
We adapt two basic collection strategies in Section~\ref{S-coll}.
We have developed two new algorithms for the full unipotent subgroup of a split 
reductive algebraic group:
The first is a new collection strategy called \emph{collection from the outside}
(Section~\ref{S-collout}). 
The second is a direct method for computing products and inverses in 
classical groups using standard representations 
(Section~\ref{S-class}).
Section~\ref{S-analysis} gives the asymptotic analysis and  the proof of
Theorem~\ref{T-main}.  
Section~\ref{S-time} compares practical timings for the various methods
considered and describes the default method used in Magma \cite{BosmaCannonPlayoust97} from
version~2.13.

\section{Unipotent groups and presentations}\label{S-pres}
Throughout this paper, $\F$ is a field and $\E$ is a commutative unital 
algebra over $\F$.
We assume that we have effective algorithms for the basic element operations in
$\F$ and $\E$.
We find it convenient to use the scheme-theoretic definition of algebraic
groups.
So an algebraic group defined over $\F$ is a functor from the category of all
commutative unital algebras over $\F$ to the category of groups, satisfying the
appropriate additional conditions \cite{DemazureGabriel70Book,Waterhouse79Book}.
If $G$ is an algebraic group, then $G(\E)$ is an abstract group, called the
\emph{rational point group} of $G$ over $\E$.
For example, the additive group $\Ga$ has rational point group $\Ga(\E)=\E^+$.

We define an \emph{$\F$-unipotent group}
to be an algebraic group defined over $\F$ with a normal series in
which every quotient is isomorphic to $\Ga$.
This is analogous to the definition of an $\F$-soluble algebraic group in 
\cite{DemazureGabriel70Book}.
\begin{Prop}\label{P-Funip}
Over a perfect field $\F$, every connected unipotent group is
$\F$-unipotent.
Over a field $\F$ of characteristic zero, every unipotent group is
$\F$-unipotent.
\end{Prop}
\begin{proof}
The first statement follows from 
\cite[Proposition~5, Corollary~2]{Rosenlicht57}.
The second follows from the first using the standard fact that all unipotent
groups are connected in characteristic zero.
\end{proof}
\noindent The
group $\boldsymbol{\alpha}_p$, defined in \cite{Waterhouse79Book} over fields of
characteristic $p>0$, is unipotent but not $\F$-unipotent.

Let $U$ be an $\F$-unipotent group.  
Fix a central series
$$
  U=U_1>U_2>\dots>U_{N+1}=1,
$$
such that each $U_r/U_{r+1}$ is $\F$-isomorphic to $\Ga$.
The projection $U_r\to U_r/U_{r+1}\iso\Ga$ splits as an $\F$-morphism of schemes
by~\cite[Theorem~16.2.6]{Springer98Book}.
Fix splitting maps $x_r:\Ga\to U_r$.
Clearly $U$ is parametrised by $N$-dimensional affine space:
$$
  \A^N\to U,\quad (a_1,\dots,a_N)\mapsto x_1(a_1)\cdots x_N(a_N).
$$
Multiplication and inversion in $U(\E)$ are given by polynomials.
To be precise
\begin{align}
  \prod_{r=1}^N x_r(a_r) \prod_{r=1}^N x_r(b_r) &= 
    \prod_{r=1}^N x_r(F_r(a_1,\dots,a_N,b_1,\dots,b_N)),\label{E-defmult}\\
  \left( \prod_{r=1}^N x_r(a_r) \right)^{-1} &= 
    \prod_{r=1}^N x_r(G_r(a_1,\dots,a_N)),\label{E-definv}
\end{align}
where all products are written in ascending order,
each $F_r$ is a polynomial in $2N$
indeterminates, and each $G_r$ is a polynomial in $N$ indeterminates.
The $F_r$ and $G_r$ are called \emph{Hall polynomials}  \cite{Hall69}.
These polynomials have coefficients in $\F$, but 
tend to be very large and unwieldy.
In order to do practical computations in $U$,
we need a more concise description.

We now construct a presentation for the rational point group $U(\E)$.
Since the series $U(\E)=U_1(\E)>U_2(\E)>\dots>U_N(\E)>U_{N+1}(\E)=1$ is central, we have
$$
  \begin{array}{rcr}
  x_r(a)x_r(b) &\in& x_r(a+b)U_{r+1}(\E),\\
  x_r(a)^{-1}  &\in& x_r(-a)U_{r+1}(\E),\\ 
  x_s(b)x_r(a) &\in& x_r(a)x_s(b)U_{s+1}(\E),
  \end{array}
$$
for $a,b\in\E$ and $1\le r<s\le N$.
Hence we have relations
\begin{align}
  x_r(a)x_r(b) &= x_r(a+b)    \prod_{t=r+1}^{N} x_t(f_{rt}(a,b)), \label{E-r}\\
  x_r(a)^{-1}  &= x_r(-a)     \prod_{t=r+1}^{N} x_t(g_{rt}(a)),   \label{E-inv}\\
  x_s(b)x_r(a) &= x_r(a)x_s(b)\prod_{t=s+1}^{N} x_t(h_{rst}(a,b)),\label{E-rs}
\end{align}
where $f_{rt}$, $g_{rt}$, and $h_{rst}$ are polynomials defined over $\F$.
We note that many of these relations are redundant (including \eqref{E-inv} for
all $r$), but the extra relations are useful for computation.
\begin{Thm}\label{T-FC}
Let $U$ be an\/ $\F$-unipotent group and let\/ $\E$ be a commutative $\F$-algebra.
Let $\tilde{U}(\E)$ be the group with generators 
$x_r(a)$, for $a\in\E$, $r=1,\dots,N$, and relations~(\ref{E-r}), (\ref{E-inv}),
and~(\ref{E-rs}), for $a,b\in\E$, $1\le r<s\le N$.
Then the natural map $\tilde{U}(\E)\to U(\E)$ is an isomorphism of abstract groups.
\end{Thm}
\begin{proof}
Since the given relations hold in $U(\E)$, the map is well defined.
The map is onto because $U(\E)$ is generated by the images of the generators of 
$\tilde{U}(\E)$.
Every element of $\tilde{U}(\E)$ is a word with terms of the form $x_r(a)$ or
$x_r(a)^{-1}$.
This word can be \emph{collected} into a product $\prod_{r=1}^N x_r(a_r)$.
This is achieved by first eliminating all
inverses using (\ref{E-inv}), then putting the terms in order by the subscripts
using (\ref{E-rs}) and removing multiple terms with the same subscript using 
(\ref{E-r}).
If the words $\prod_{r=1}^N x_r(a_r)$ and $\prod_{r=1}^N x_r(b_r)$ are equal in 
$U(\E)$, then $a_r=b_r$ for all $r$, and so these words are also equal in
$\tilde{U}(\E)$.
Hence the map is injective and we are done.
\end{proof}
\noindent We say that the group scheme $U$ is \emph{presented by} $\tilde{U}$.

Now suppose we are given an arbitrary system of $\F$-polynomials
$f_{rt}(a,b)$ and %, for $1\le r<t\le N$;
$g_{rt}(a,b)$, for $1\le r<t\le N$; and 
$h_{rst}(a)$, for $1\le r<s<t\le N$.
Define the group functor $\tilde{U}$ by taking
$\tilde{U}(\E)$ to be the abstract group given by generators 
$x_r(a)$, for $a\in\E$, $r=1,\dots,N$, and relations~(\ref{E-r}), (\ref{E-inv}), 
and~(\ref{E-rs}).
We call $\tilde{U}$ an \emph{FC group functor over $\F$}.
FC stands for field-commutator, since the relations involve field operations
and commutators, just as the PC presentation of a nilpotent group involves
powers and commutators.
We call $\tilde{U}$ \emph{consistent} if the map 
$$
  \E^N\to\tilde{U}(\E), \quad (a_1,\dots,a_N)\mapsto x_1(a_1)\cdots x_N(a_N)
$$
is injective for every commutative $\F$-algebra $\E$.
Theorem~\ref{T-FC} implies that every $\F$-unipotent group is presented by a
consistent FC group functor.
We now prove the converse:
\begin{Thm}\label{T-FCconv}
Every consistent FC group functor defined over $\F$ is an $\F$-unipotent group.
\end{Thm}
\begin{proof}
Let $\tilde{U}$ be the consistent FC group functor.
Ignoring the multiplication, we can consider $\tilde{U}$ to be the 
$N$-dimensional affine scheme.
Using collection, as in the proof of Theorem~\ref{T-FC}, we can find polynomials
$F_r$ and $G_r$ such that equations~(\ref{E-defmult}) and~(\ref{E-definv}) are
satisfied in $\tilde{U}(\E)$.
So $\tilde{U}$ is an $\F$-algebraic group scheme,
since $F_r$ and $G_r$ are clearly defined over $\F$.
Finally define algebraic subgroups $U_r = \prod_{k=r}^N \im(x_r)$.
These give a normal series for $U$ in which
every quotient is isomorphic to $\Ga$, and so $\tilde{U}$ is an $\F$-unipotent group. 
\end{proof}

Let $U$ be an $\F$-unipotent group.
Suppose the projection $U_r\to\Ga$ splits as a homomorphism of $\F$-group
schemes, not just as a morphism of $\F$-schemes.
Then we can take $x_r:\Ga\to U_r$ to be a homomorphism, and so replace 
(\ref{E-r}) and~(\ref{E-inv}) by
\begin{align}
  x_r(a)x_r(b) &= x_r(a+b). \label{E-splitr}
\end{align}
It follows immediately that
\begin{align}
  x_r(a)^{-1}  &= x_r(-a), \label{E-splitinv}
\end{align}
and so all the polynomials $f_{rt}$ and $g_{rt}$ are zero.
If $x_r$ is a homomorphism for $r=1,\dots,N$, we call the corresponding
FC group functor \emph{split}.
\begin{Thm}
If\/ $\F$ is a field of characteristic zero, then every unipotent group defined
over $\F$ is
presented by a split FC group functor over\/ $\F$.
\end{Thm}
\begin{proof}
Since $\F$ has characteristic zero, every unipotent group $U$ defined over $\F$ is 
$\F$-unipotent by Proposition~\ref{P-Funip}.
By induction, we can assume that $U/U_N$ is presented by a split FC group
functor.
Fix maps $y_r:\Ga\to U/U_N$, for $r=1,\dots,N-1$, defining this functor.
By \cite[Proposition~VII.8]{Serre88Book}, $\Ext(\Ga,\Ga)=0$.
Hence $\Ext(U/U_N,\Ga)=0$, by repeated application of the long exact sequence
for $\Hom(\circ,\Ga)$.
So there exists a homomorphism $y:U/U_N \to U$ splitting the projection
$U\to U/U_N$.
Define $x_r=y\circ y_r$ for $r=1,\dots,N-1$.
Take $x_N$ to be the $\F$-injection $\Ga\iso U_N\to U$.
These maps clearly define a split FC group functor presenting $U$.
\end{proof}
\noindent If $\F$ has positive characteristic, then the Witt-vector groups 
\cite{Serre88Book}
provide examples of $\F$-unipotent groups which cannot be presented by
split FC group functors.
The full unipotent subgroup of a split reductive group is always presented by 
a split FC group functor, as we show in Proposition~\ref{P-redsplit} below.

\section{Collection and symbolic collection}\label{S-coll}
We now extend some of the standard collection strategies for PC groups to 
FC group functors.
The precise order in which the relations are applied has a huge impact on the
speed of collection.
Many strategies have been suggested, and we have not attempted to extend them
all to FC group functors.
We have implemented two fundamental techniques:
\emph{collection from the left}
\cite{Leedham-GreenSoicher90,Vaughan-Lee90}; and a slightly improved version of
\emph{collection to the left} \cite{Hall-Jr76Book} 
(we collect the rightmost rather than the leftmost occurrence of the
least uncollected letter).

Let $U$ be an FC group functor over $\F$, and let $\E$ be a commutative
$\F$-algebra.
The algorithms in this section operate on a word $w\in U(\E)$.
This word is always equal to 
$\prod_{i=1}^{M} x_{r_i}(a_i)^{\ep_i}$,
ie, the parameters $M\in\N$, $a_i\in\E$, $\ep_i=\pm1$, and $r_i\in\{1,\dots,N\}$ 
are automatically modified when $w$ is.
Algorithm~\ref{A-collsub} describes the basic step of collection.
When we say ``\emph{apply} a certain relation to a subword'', we mean match the 
subword with the left hand side of the relation, and replace it by the 
right hand side. 
{\sc CollectSubword} looks at the term at position $j$ in the word $w$, and 
either removes an inverse (if $\ep_j=-1$) or ensures that $r_{j-1}<r_j$.
In addition to the modified word $w$, it returns  indices $j_1$ and $j_2$.
\begin{algorithm}
\begin{algorithmic}
\STATE{\hspace{-1em}$\text{\sc CollectSubword} := \text{\bf function}
  (U, w=\prod_{i=1}^{M} x_{r_i}(a_i)^{\ep_i}, j)$}\\
      \IF{$\ep_j=-1$}
        \STATE{apply \eqref{E-inv} to the subword $x_{r_j}(a_j)^{-1}$} 
%       \STATE{$j_1:=j,\quad j_2:=j_1+\#\{t:g_{r_jt}(a_j)\ne0\}$}
        \STATE{$j_1:=j,\quad j_2:=j$}
      \ELSIF{$j>1$ {\bf and} $r_{j-1}=r_j$}
	\STATE{apply \eqref{E-r} to the subword $x_{r_{j}}(a_{j-1}) x_{r_j}(a_j)$}
	\STATE{$j_1:=j-1,\quad j_2:=j_1+\#\{t:f_{r_jt}(a_{j-1},a_j)\ne0\}$}
      \ELSIF{$j>1$ {\bf and} $r_{j-1}>r_j$} %\COMMENT{apply \eqref{E-rs}}
	\STATE{apply \eqref{E-rs} to the subword 
	  $x_{r_{j-1}}(a_{j-1})	x_{r_j}(a_j)$}
	\STATE{$j_1:=j-1$}
        \IF{$j_1>1$ {\bf and} $r_{j_1-1} < r_{j_1}$}
          \STATE{$j_2:=j_1$}
        \ELSE
          \STATE{$j_2:=j_1+1+\#\{t:h_{r_jr_{j-1}t}(a_{j-1},a_j)\ne0\}$}
        \ENDIF
      \ELSE
        \STATE{$j_1:=j,\quad j_2:=j_1+1$}
      \ENDIF
  \RETURN $w$, $j_1$, $j_2$\\
\end{algorithmic}
\caption{Collect subword}
\label{A-collsub}
\end{algorithm}
Collection to the left (Algorithm~\ref{A-colltoleft}) works by collecting all
terms $x_1(a)$, followed by all terms $x_2(a)$, and so on.
This uses the index $j_1$, which gives the new largest $j$ such that $r_j=r$.
\begin{algorithm}
\begin{algorithmic}
\REQUIRE An FC group functor $U$ and a word $w=\prod_{i=1}^{M} x_{r_i}(a_i)^{\ep_i}$.
\ENSURE  A product $\prod_{r=1}^{N} x_{r}(b_r)$
  that is equal to $w$ as an element of $U(\E)$.
  \FOR{$r := 1$ to $N$}
    \STATE{let $j$ be the largest $i$ such that $r_i=r$}\\
    \WHILE{$j \ge r$}
      \STATE{$w,j_1,j_2 := \text{\sc CollectSubword}(U,w,j),\quad j:=j_1$}\\
    \ENDWHILE
  \ENDFOR
  \RETURN $w$
\end{algorithmic}
\caption{Collection to the left}
\label{A-colltoleft}
\end{algorithm}
Collection from the left (Algorithm~\ref{A-collfromleft}) goes through the word from left to right, correcting
each term that is out of position.
This uses the index $j_2$, which gives the next term which is potentially out of
position.
\begin{algorithm}
\begin{algorithmic}
\REQUIRE An FC group functor $U$ and a word $w=\prod_{i=1}^{M} x_{r_i}(a_i)^{\ep_i}$.
\ENSURE  A product $\prod_{r=1}^{N} x_{r}(b_r)$
  that is equal to $w$ as an element of $U(\E)$.
  \STATE{$j := 1$}\\
  \WHILE{$j\le M$}
    \STATE{$w,j_1,j_2 := \text{\sc CollectSubword}(U,w,j),\quad j:=j_2$}\\
  \ENDWHILE
  \RETURN $w$
\end{algorithmic}
\caption{Collection from the left}
\label{A-collfromleft}
\end{algorithm}

\emph{Symbolic collection} is a standard method for improving the efficiency of element
multiplication in a PC group.
This depends on the observation that we can collect a generic product,
and then substitute into polynomials for subsequent collections.
This is particularly easy in the case of FC group functors:
simply take $\E=\F[a_1,a_2,\dots,a_N,b]$ and do $N$ collections in $U(\E)$ 
to get relations
\begin{align}
  \bigg(\prod_{s=r+1}^N x_s(a_s)\bigg)x_r(b) &= 
  x_r(b)\bigg(\prod_{s=r+1}^N x_s(c_{rs})\bigg) \label{E-symb}
\end{align}
where $c_{rs}$ is a polynomial in $b$ and $a_{r+1},\dots,a_N$.
Now, taking arbitrary $\E$ again,
we can multiply two collected words
$\prod_{i=1}^N x_i(a_i)$ and $\prod_{i=1}^N x_i(b_i)$
by substituting values from $\E$ into the $(N-1)N/2$ polynomials $c_{rs}$ 
in the obvious manner.
A similar method can be used to compute inverses.

The advantage of symbolic collection is that each operation is faster.
The disadvantage is that more preprocessing time and memory are required.
In order to save memory, we represent our polynomials as straight-line programs
(see \cite{Leedham-Green01} for a description of straight-line programs for group
elements;
the implementation for polynomials is due to Allan Steel).
This means that the  polynomials are basically just the collection preserved
in amber, so the collection method used is still important. 
We note that there is another common symbolic collection algorithm, 
called Deep Thought \cite{Leedham-GreenSoicher98,Merkwitz97},
but we have not implemented it for unipotent groups.

\section{Collection in the full unipotent subgroup}\label{S-collout}
We now describe a new collection strategy for the full unipotent subgroup of
a reductive group.
Let $G$ be an $\F$-split reductive algebraic group \cite{Springer98Book}.
Fix a split maximal torus $T$ in $G$, and a Borel subgroup $B$ containing $T$.
Let $U$ be the unipotent radical of $B$.
Since $U$ is unique up to $G$-conjugacy, we refer to $U$ as the
\emph{full unipotent subgroup} of $G$.
Let $\Phi$ be the root system of $G$ with respect to $T$,
and let $\Phi^+\subseteq\Phi$ be the positive roots with respect to $B$.

Write $\Phi^+=\{\al_1,\al_2,\dots,\al_N\}$ with the roots in an \emph{order
compatible with height}, ie, $\height(\al_r)<\height(\al_s)$ implies
$r<s$. 
For each $\al\in\Phi^+$, there is a subgroup $X_\al$ of $U$
isomorphic to $\Ga$.
Write $x_r$ for the isomorphism $\Ga\to X_{\al_r}$.
Then $U_r=\prod_{s=r}^N X_{\al_s}$, for $r=1,\dots,N+1$, defines a central series for $U$.
The corresponding FC group functor has relations (\ref{E-splitr}) and
\begin{align}
  x_s(b) x_r(a)  &= x_r(a) x_s(b) 
    \prod_{\substack{\al_t=i\al_r+j\al_s,\\i,j>0}}
     x_t(C_{ij\al_r\al_s}a^ib^j), \label{E-redrs}
\end{align}
where the constants $C_{ij\al_r\al_s}$ are defined as in \cite{Carter72Book}.
Recall that these constants depend on the combinatorics of $\Phi$, and on the 
choice of a sign for each nonsimple positive root.
%These \emph{extraspecial signs} are used in Section~\ref{S-class} below.
In \cite{CohenMurrayTaylor04}, a method for computing these constants is given 
which is efficient for small rank groups.
For large ranks, we outline a new method in Section~\ref{S-class}.
We now have:
\begin{Prop}\label{P-redsplit}
The full unipotent subgroup of a split reductive group is presented by a split FC
group functor.
\end{Prop}

Note that most of the polynomials $h_{rst}$ from (\ref{E-rs}) of
Section~\ref{S-pres} are zero in this case.
This gives us much greater flexibility in how we collect words.
The ordering of the roots is of vital importance in this section.
We specify an ordering in terms of subscripts: $\al_1,\dots,\al_N$.
The subscripts on the injections $x_r:\Ga\to U$ are always kept in agreement
with the root ordering under discussion.

Words in $U$ need not be collected into an order compatible with
height.
In fact, the algorithms of the previous section work for all orderings
$\al_1,\al_2,\dots,\al_N$ of $\Phi^+$ with the property that
$\al_r+\al_s=\al_t$ implies $t>r$ and $t>s$.
We call such an ordering \emph{left-additive}.

An ordering $\al_1,\al_2,\dots,\al_N$ of $\Phi^+$ is called 
\emph{additive} if $\al_r+\al_s=\al_t$ implies $t$ lies between $r$ and $s$ (ie,
$r<t<s$ or $s<t<r$). 
Additive orderings reflect more of the combinatorial structure of the
root system than left-additive orderings do.
The existence and construction of additive orderings is considered in
\cite{Papi94} (see Section~\ref{S-analysis} below for more details).

In order to collect a word into the additive ordering, we replace relation
\eqref{E-redrs} with
\begin{align}
  x_s(b) x_r(a) %&= x_r(a)\big(x_r(a)^{-1}x_s(-b)^{-1}x_r(a)x_s(-b)\big)x_s(b)\\
  &=  x_r(a) \bigg( 
        \prod_{\substack{\al_t=i\al_r+j\al_s,\\i,j>0}}
            x_t(C_{ji\al_s\al_r}a^i(-b)^j) 
    \bigg) 
    \;x_s(b), \label{E-addrs}
\end{align}
which can be proved by applying \eqref{E-redrs} to $x_r(a)x_s(-b)$.
We then use \emph{collection from the outside} (Algorithm~\ref{A-collout}),
which is a modified version of collection from the left.
The basic idea is to run collection %from the left 
from both sides simultaneously, until we meet in the middle.
The two collections could easily be run in parallel, but we alternate between
them. The subroutine {\sc CollectSubword} from Algorithm~\ref{A-collfromleft}
is slightly modified for both collections ({\sc CollectSubwordL} and {\sc CollectSubwordR}). 
The returned value $L$ is the increase in the word length and $k$ is the index of the 
next term potentially out of position.

\begin{algorithm}
\begin{algorithmic}
\STATE{\hspace{-1em}$\text{\sc CollectSubwordL} := \text{\bf function}
  (U,w=\prod_{i=1}^{M} x_{r_i}(a_i)^{\ep_i}, j)$}\\
      \IF{$\ep_j=-1$}
        \STATE{apply \eqref{E-splitinv} to the subword $x_{r_j}(a_j)^{-1}$} \\
	\STATE{$k:=j, \quad L:=0$}
      \ELSIF{$j>0$ {\bf and} $r_{j-1}=r_j$}
	\STATE{apply \eqref{E-splitr} to the subword $x_{r_{j}}(a_{j-1}) x_{r_j}(a_j)$}\\
	\STATE{$k:=j-1,\quad L:=-1$}
      \ELSIF{$j>0$ {\bf and} $r_{j-1}>r_j$}
	\STATE{apply \eqref{E-addrs} to the subword 
	  $x_{r_{j-1}}(a_{j-1})	x_{r_j}(a_j)$}\\
	\STATE{$k:=j-1,\quad 
	  L:=\#\{t:\al_t=k\al_{r_j}+l\al_{r_{j-1}}\text{ for }k,l>0\}$}
      \ELSE
        \STATE{$k:=j+1,\quad L:=0$}
      \ENDIF
  \RETURN  $w$, $k$, $L$
\end{algorithmic}
\vspace{0.5cm}
\begin{algorithmic}
\STATE{\hspace{-1em}$\text{\sc CollectSubwordR} := \text{\bf function}
  (U,w=\prod_{i=1}^{M} x_{r_i}(a_i)^{\ep_i}, j)$}\\
      \IF{$\ep_j=-1$}
        \STATE{apply \eqref{E-splitinv} to the subword $x_{r_j}(a_j)^{-1}$} \\
	\STATE{$k:=j, \quad L:=0$}
      \ELSIF{$j<M$ {\bf and} $r_j=r_{j+1}$}
	\STATE{apply \eqref{E-splitr} to the subword $x_{r_{j}}(a_{j}) x_{r_j}(a_{j+1})$}\\
	\STATE{$k:=j,\quad L:=-1$}
      \ELSIF{$j<M$ {\bf and} $r_{j}>r_{j+1}$}
	\STATE{apply \eqref{E-addrs} to the subword 
	  $x_{r_j}(a_j)   x_{r_{j+1}}(a_{j+1}) $}\\
	\STATE{$L:=\#\{t:\al_t=k\al_{r_{j+1}}+l\al_{r_j}\text{ for }k,l>0\}$,
            \quad $k := j+1+L$}
      \ELSE
        \STATE{$k:=j-1,\quad L:=0$}
      \ENDIF
  \RETURN  $w$, $k$, $L$
\end{algorithmic}
\vspace{0.5cm}
\begin{algorithmic}
\REQUIRE An FC group functor $U$ and a word $w=\prod_{i=1}^{M} x_{r_i}(a_i)^{\ep_i}$.
\ENSURE  A product $\prod_{r=1}^{m} x_{r}(b_r)$
  that is equal to the input as an element of $U(\E)$.
  \STATE{$i:=1,\quad j := M$}\\
  \WHILE{$i<j$}
    \STATE{$w,k,L := \text{\sc CollectSubwordL}(U,w,i),\quad   i:=k,\  M:= M+L,\ j:= j+L $}\\
    \IF{$i<j$}
      \STATE{$w,k,L := \text{\sc CollectSubwordR}(U,w,j),\quad j:=k,\  M:= M+L $}
    \ENDIF
  \ENDWHILE
  \RETURN  $w$
\end{algorithmic}
\caption{Collection from the outside}
\label{A-collout}
\end{algorithm}

Finally we note that symbolic collection works with collection from the outside, 
with the obvious minor modifications.

\section{The full unipotent subgroup of a classical group}\label{S-class}
In this section, we present an alternative to collection which is much more
efficient when $G$ has large semisimple rank. 
We derive formulas for the defining polynomials of the
full unipotent subgroup $U$ of a split classical group $G$.
This allows us to write code that implicitly applies the defining polynomials of
$U$ without having to  store them explicitly in memory.
We do not consider exceptional groups here, because their semisimple rank
is at most $8$.

The rough outline of our method is as follows:
We index the roots by pairs of integers.
We then construct a minimal degree matrix representation of $G$.
We take the basis for the representation consisting of weight vectors, ordered
according to the dominance ordering  on the corresponding weights
\cite{Humphreys78Book}.
We then order the roots by going down each column of this matrix
representation, and seeing where the parameters corresponding to each root
appear.
We call this the \emph{representation ordering}.
Note that the representation ordering is also an additive ordering,
except in Cartan type $\CC_\ell$, where it gives an additive ordering on the
coroots. 

The representation of roots by pairs of integers can also be used to compute
the constants $C_{ij\al_r\be_r}$. 
For large classical groups, this is much more efficient than the method given in
\cite{CohenMurrayTaylor04}.
We omit the details since this is a technical but straightforward application of
formulas in \cite{Carter72Book}.

We consider each classical type in Subsections~\ref{SS-A} to~\ref{SS-D}.
In each subsection, we fix a particular isogeny type and particular extraspecial
signs \cite{Springer98Book}.
The choice of isogeny class is irrelevant, since the structure of the unipotent
subgroup is independent of isogeny.
We can easily transform between different choices of extraspecial signs
using Theorem~29 of~\cite{Steinberg68}.

\subsection{Cartan type $\CA_\ell$: Linear of degree $\ell+1$}\label{SS-A}
Let $G=\SL_{l+1}$ and let $U$ be the algebraic subgroup of all lower
unitriangular matrices.
The root system of $G$ has Cartan type $\CA_\ell$.
Let $V=\R^{\ell+1}$ with basis $e_1,\dots,e_{\ell+1}$.
Then the roots are
$$
  \al_{ij}=e_i-e_j
$$
for $i,j=1,\dots,\ell+1$ with $i\ne j$.
The simple roots are $\al_{i,i+1}$ for $i=1,\dots,\ell$.
A root $\al_{ij}$ is positive if, and only if, $i<j$.
Roots add by the formulas $\al_{ij}+\al_{jm}=\al_{im}$, and
$\al_{ij}+\al_{km}=0$ unless $j=k$ or $i=m$.

Define the map $x_{ij}:\Ga\to G$ by $x_{ij}(a)=I+aE_{ji}$.
The representation order on the roots is simply the lexicographic order on the
corresponding pairs of integers.
We label the coordinates of $\A^{(\ell+1)\ell/2}$ by the root pairs in
this order, ie, $\bfa\in \A^{(\ell+1)\ell/2}(\E)$ has the form
$$
  \bfa=(a_{12},a_{13},\dots,a_{1,\ell+1},\; a_{23},\dots,a_{2,\ell+1},\; \dots,\; 
    a_{\ell,\ell+1}).
$$
We can now define a parametrisation $\varphi:\A^{(\ell+1)\ell/2}\to U$ by
$$
  \varphi(\bfa) := \prod_{i=1}^{\ell+1}\prod_{j=i+1}^{\ell+1} x_{ij}(a_{ij}) =
  \left(\begin{array}{ccccc}
    1\\
    a_{12}	 &1\\
    a_{13}	 & a_{23}	&1\\
    \vdots	 & \vdots	&\ddots &\ddots\\
    a_{1,\ell+1} & a_{2,\ell+1} &\dots  & a_{\ell,\ell+1} &1
  \end{array}\right).
$$
We now get $\varphi(\bfa)\varphi(\bfb)=\varphi(\bfc)$ where
$$
  c_{ij} = a_{ij} + \sum_{i<k<j}b_{ik}a_{kj} + b_{ij}.
$$
Also $\varphi(\bfa)^{-1}=\varphi(\bfd)$ where
$$
  d_{ij} = -a_{ij} -\sum_{i<k<j}d_{ik}a_{kj}
$$
The formulas for inversion are defined recursively and are computed in
reverse representation order.
We note that it is easy to derive a direct formula for $d_{ij}$,
but the recursive version can be evaluated with fewer operations.

\subsection{Cartan type $\CB_\ell$: Orthogonal of degree $2\ell+1$}\label{SS-B}
Let $F_m$ be the $m\times m$ matrix over $\F$ of the form
$$
  \left(\begin{matrix}
    0&\dots&0&1\\
    0&\dots&1&0\\
    \vdots&\addots&\vdots&\vdots\\
    1&\dots&0&0
  \end{matrix}\right).
$$
We assume, for this subsection only, that $\F$ has odd characteristic.  
Since the group of type $\CB_\ell$ is isomorphic to the group of type 
$\CC_\ell$ in characteristic $2$, this restriction is not critical.
Let $G=\SO_{2\ell+1}$ be the special orthogonal group of the orthogonal form
with matrix 
$$
  \left(\begin{matrix}
    0&0&F_\ell \\ 0&2&0 \\ F_\ell&0&0
  \end{matrix}\right).
$$
Let $U$ be the group of all lower unitriangular matrices in $G$.
The root system of $G$ has Cartan type $\CB_\ell$.
Let $V=\R^{\ell}$ with basis $e_1,\dots,e_{\ell}$.
The roots are
$$
  \al_{si,tj}= se_i-te_j\quad\text{and}\quad\al_{si,0}=se_i,
$$
for $i,j=1,\dots,\ell$ with $i\ne j$ and $s,t=\pm1$.
For the sake of readability, we write $\ibar$ instead of $-i$ in subscripts,
eg, $\al_{2,-3}$ is denoted $\al_{2\m3}$.
Note that $\al_{ij}=\al_{\jbar\ibar}$ for all $i,j=\pm1,\dots,\pm\ell$.
The simple roots are $\al_{i,i+1}$, for $i=1,\dots,\ell-1$, 
and $\al_{n0}$.
A root $\al_{i,tj}$ for $i,j>0$ is positive if, and only if, $i<j$;
a root $\al_{si,0}$ is positive if, and only if, $s=+1$.

Define the root maps by
\begin{align*}
  x_{ij}(a) &= I+a(E_{ji} - E_{2\ell-i+2,2\ell-j+2}),\\
  x_{i\jbar}(a) &= I+a(E_{2\ell-j+2,i} - E_{2\ell-i+2,j}),\\
  x_{i0}(a)     &= I+a(2E_{\ell+1,i} - E_{i,\ell+1}) 
    - a^2E_{2\ell-i+2,i},\quad\text{and}\\
  x_{\ibar0}(a) &= I+a(E_{2\ell-i+2,\ell+1}-2E_{\ell+1,2\ell-i+2}) 
    - a^2E_{i,2\ell-i+2}.
\end{align*}
Let $J_i$ be the sequence of integers
$[i+1,i+2,\dots,\ell,0,-\ell,\dots,-(i+2),-(i+1)]$.
Let $J'_{i} := J_i \setminus \{0\}$.
The representation order on the positive roots is the lexicographic order on 
pairs, with the integers ordered as in $J_0$.
Label the coordinates of $\A^{\ell^2}$ by the root pairs in
this order.
Now define a parametrisation $\varphi:\A^{\ell^2}\to U$ by
\begin{align*}
  \varphi(\bfa) &:= \prod_{i=1}^{\ell}\prod_{j\in J_i} x_{ij}(a_{ij})\\
  &=\left(\begin{array}{ccccccccc}
    1\\
    a_{12}	& 1\\
    \vdots	& \ddots    & \ddots\\
    a_{1\ell}	& \dots     & a_{\ell-1,\ell}	& 1\\
    2a_{10}	& \dots     & 2a_{\ell-1,0}	& 2a_{\ell0}	     & 1\\
    a_{1\m\ell} & \dots     & a_{\ell-1,\m\ell} & a''_\ell	     & a'_{\ell0}    & 1 \\
    \vdots	& \addots   & \addots           & a'_{\ell-1,\m\ell} & a'_{\ell-1,0} & a_{\ell-1,\ell} & \ddots \\
    a_{1\m2}	& a''_2     & \addots           & \vdots	     & \vdots        & \vdots          & \ddots & 1 \\
    a''_1       & a'_{1\m2} & \dots             & a'_{1\m\ell}	     & a'_{10}       & a'_{1\ell}      & \dots  & a'_{12} & 1 \\
  \end{array}\right).
\end{align*}
where
\begin{align*}
  a''_{i}     &= -2{a_{i0}}^2 - \sum_{i<k\le\ell} a_{ik}a_{i\m{k}},\qquad
  a'_{i0}     = -a_{i0} -\sum_{i<k\le\ell} a_{ik}a'_{k0}, \\
  a'_{ij}     &= -a_{ij} -\sum_{i<k<j} a_{ik}a'_{kj}, \\
  a'_{i\jbar} &= -a_{i\jbar} 
        - \sum_{i<k<j} a_{ik}a'_{k\jbar}
        - a_{ij}a''_j
        - 2a_{i0}a_{j0}
        - \sum_{k\in J'_j} a_{ik}a_{j\kbar} .
\end{align*}
Now $\varphi(\bfa)\varphi(\bfb)=\varphi(\bfc)$ where
\begin{align*}
  c_{i0} &= a_{i0} + b_{i0} +\sum_{i<k\le\ell}b_{ik}a_{k0} ,  \qquad
  c_{ij} = a_{ij} + b_{ij} +\sum_{i<k<j}     b_{ik}a_{kj} ,  \\
  c_{i\jbar} &= a_{i\jbar} + b_{i\jbar} 
                            + \sum_{i<k<j}     b_{ik}a_{k\jbar}
                            + a''_{j}b_{ij}
                            + 2b_{i0}a'_{j0}
                            + \sum_{k\in J'_j} b_{ik}a'_{j\kbar}.
\end{align*}
And $\varphi(\bfa)^{-1}=\varphi(\bfd)$ where
\begin{align*}
  d_{i0}     &= -a_{i0}     - \sum_{i<k\le\ell}d_{ik}a_{k0},\qquad
  d_{ij}     = -a_{ij}     - \sum_{i<k<j}     d_{ik}a_{kj},\\
  d_{i\jbar} &= -a_{i\jbar} - \sum_{i<k<j}     d_{ik}a_{k\jbar}
                            - a''_{j}d_{ij}
                            - 2d_{i0}a'_{j0}
                            - \sum_{k\in J'_j} d_{ik}a'_{j\kbar}.
\end{align*}
All of these equations are recursive and are computed in reverse
representation order.

\subsection{Cartan type C: Symplectic}\label{SS-C}
Let $G=\Sp_{2\ell}$ be the symplectic group of the symplectic form with matrix
$$
  \left(\begin{matrix} 0 & F_\ell \\ -F_\ell & 0 \end{matrix}\right)
$$
Let $U$ be the group of all lower unitriangular matrices in $G$.
The root system of $G$ has Cartan type $\CC_\ell$.
Let $V=\R^{\ell}$ with basis $e_1,\dots,e_{\ell}$.
The roots are
$$  
  \al_{si,tj}= se_i-te_j,
$$
for $i,j=1,\dots,\ell$ with $i\ne j$, and $s,t=\pm1$.
Once again $\al_{ij}=\al_{\jbar\ibar}$ for all $i,j=0,\pm1,\dots,\pm\ell$.
The simple roots are $\al_{i,i+1}$ for $i=1,\dots,\ell-1$, 
and $\al_{\ell\bar{\ell}}$.
A root $\al_{i,tj}$ is positive if, and only if, $i<j$.

Define the root maps by
\begin{align*}
  x_{ij}(a) &= I+a(E_{ji} - E_{2\ell-i+1,2\ell-j+1}),\\
  x_{i\jbar}(a) &= I+a(E_{2\ell-j+1,i} + E_{2\ell-i+1,j}),\quad\text{and}\\
  x_{i\ibar}(a) &= I+aE_{2\ell-i+1,i}.
\end{align*}
The representation order on the positive roots is the lexicographic order on 
pairs, with the integers ordered as in $J'_0$.
Label the coordinates of $\A^{\ell^2}$ by the root pairs in
this order.
Now define a parametrisation $\varphi:\A^{\ell^2}\to U$ by
\begin{align*}
  \varphi(\bfa) &:= \prod_{i=1}^{\ell+1}\prod_{j\in J'_i} x_{ij}(a_{ij})\\
  &=\left(\begin{array}{ccccccccc}
    1\\
    a_{12}      & 1\\
    \vdots      & \ddots     & \ddots\\
    a_{1\ell}   & \dots      & a_{\ell-1,\ell}   & 1\\
    a_{1\m\ell} & \dots      & a_{\ell-1,\m\ell} & a''_{\ell}         & 1 \\
    \vdots      & \addots    & \addots  	 & a'_{\ell-1,\m\ell} & a'_{\ell-1,\ell} & \ddots \\
    a_{1\m2}    & a''_{2}    & \addots           & \vdots	      & \vdots  	 & \ddots & 1 \\
    a''_{1}     & a'_{1\m2}  & \dots		 & a'_{1\m\ell}       & a'_{1\ell}	 & \dots  & a'_{12} & 1 \\
  \end{array}\right),
\end{align*}
where
\begin{align*}
  a''_{i}     &= a_{i\ibar} - \sum_{i<j\le\ell} a_{ij}a_{i\jbar}, \qquad
  a'_{ij}     = a_{ij} -\sum_{i<k<j} a_{ik}a'_{kj},\\
  a'_{i\jbar} &= a_{i\jbar} - \sum_{i<k<j} a_{ik}a'_{k\jbar}  - a_{ij}a''_j
    - \sum_{k\in J'_j} \sign(k) a_{ik}a_{j\kbar} .
\end{align*}
Now $\varphi(\bfa)\varphi(\bfb)=\varphi(\bfc)$ where
\begin{align*}
  c_{ij}     &= a_{ij}     + b_{ij}     + \sum_{i<k<j}b_{ik}a_{kj},\qquad
  c_{i\ibar} = a''_{i}    + b''_{i}    + \sum_{i<k\le\ell} c_{ik}c_{i\kbar}
                                        + \sum_{k\in J'_i} b_{ik}a'_{i\kbar},\\
  c_{i\jbar} &= a_{i\jbar} + b_{i\jbar} + \sum_{i<k<j}b_{ik}a_{k\jbar}
                                        + a''_{j}b_{ij}
                                        + \sum_{k\in J'_j} b_{ik}a'_{j\kbar}.
\end{align*}
And $\varphi(\bfa)^{-1}=\varphi(\bfd)$ where
\begin{align*}
  d_{ij}     &= -a_{ij}     - \sum_{i<k<j}d_{ik}a_{kj},\qquad
  d_{i\ibar} = -a''_{i}    + \sum_{i<k\le\ell}  d_{ik}d _{i\kbar}
                            - \sum_{k\in J'_{i}}d_{ik}a'_{i\kbar},\\
  d_{i\jbar} &= -a_{i\jbar} - \sum_{i<k<j}d_{ik}a_{k\jbar}
                            - a''_{j}d_{ij}
                            - \sum_{k\in J'_{j}}d_{ik}a'_{j\kbar}. \\
\end{align*}

\subsection{Cartan type D: Even-degree orthogonal}\label{SS-D}
Let $G=\SO_{2\ell}$ be the orthogonal group of the orthogonal form with matrix
$F_{2\ell}$.
Let $U$ be the group of all lower unitriangular matrices in $G$.
The root system of $G$ has Cartan type $\CD_\ell$.
Let $V=\R^{\ell}$ with basis $e_1,\dots,e_{\ell}$.
The roots are
$$
  \al_{si,tj}= se_i-te_j,
$$
for $i,j=1,\dots,\ell$ with $i\ne j$ and $s,t=\pm1$.
The simple roots are $\al_{i,i+1}$, for $i=1,\dots,\ell-1$, 
and $\al_{n-1,\bar\ell}$.
A root $\al_{i,tj}$ for $i,j>0$ is positive if, and only if, $i<j$.

Define the root maps by
\begin{align*}
  x_{ij}(a) &= I+a(E_{ji} - E_{2\ell-i+1,2\ell-j+1})\quad\text{and}\\
  x_{i\jbar}(a) &= I+a(E_{2\ell-j+1,i} - E_{2\ell-i+1,j}).
\end{align*}
The representation order on the positive roots is the lexicographic order on 
pairs, with the integers ordered as in $J'_0$.
Label the coordinates of $\A^{\ell(\ell-1)/2}$ by the root pairs in
this order.
Now define a parametrisation $\varphi:\A^{\ell(\ell-1)/2}\to U$ by
\begin{align*}
  \varphi(\bfa) &:= \prod_{i=1}^{\ell}\prod_{j\in J''_i} x_{ij}(a_{ij})\\
  &=\left(\begin{array}{ccccccccc}
    1\\
    a_{12}	& 1\\
    \vdots	& \ddots    & \ddots\\
    a_{1\ell}	& \dots     & a_{\ell-1,\ell}	& 1\\
    a_{1\m\ell} & \dots     & a_{\ell-1,\m\ell} & a''_{\ell}     & 1 \\
    \vdots	& \addots   & \addots           & a'_{\ell-1,\m\ell} & a'_{\ell-1,\ell} & \ddots \\
    a_{1\m2}	& a''_2     & \addots           & \vdots	     & \vdots	       & \ddots & 1 \\
    a''_1       & a'_{1\m2} & \dots             & a'_{1\m\ell}	     & a'_{1\ell}      & \dots  & a'_{12} & 1 \\
  \end{array}\right).
\end{align*}
where
\begin{align*}
  a''_{i} &= \sum_{i<j\le \ell} a_{ij}a_{i\jbar}, \qquad
     a'_{ij}  = -a_{ij} -\sum_{i<k<j} a_{ik}a'_{kj},\\
  a'_{i\jbar} &= -a_{i\jbar} - \sum_{i<k<j} a_{ik}a'_{k\jbar}- a''_{j}a_{ij}
    - \sum_{k\in J'_j} a_{ik}a_{j\kbar}  .
\end{align*}
Now $\varphi(\bfa)\varphi(\bfb)=\varphi(\bfc)$ where
\begin{align*}
  c_{ij}     &= a_{ij}     + b_{ij}     + \sum_{i<k<j}    b_{ik}a_{kj}        ,\\
  c_{i\jbar} &= a_{i\jbar} + b_{i\jbar} + \sum_{i<k<j}    b_{ik}a_{k\jbar}
              + a''_{j}b_{ij}           + \sum_{k\in J'_j}b_{ik}a'_{j\kbar}.
\end{align*}
And $\varphi(\bfa)^{-1}=\varphi(\bfd)$ where
\begin{align*}
  d_{ij} &= -a_{ij} -\sum_{i<k<j}d_{ik}a_{kj},\\
  d_{i\jbar} &= 
    -a_{i\jbar} 
    -\sum_{i<k<j}d_{ik}a_{k\jbar} 
    -a''_{j}d_{ij}
    -\sum_{k\in J'_{j}}d_{ik}a'_{j\kbar}.
\end{align*}

\section{Analysis and reductive groups}\label{S-analysis}
We can now give precise asymptotic timings for operations in reductive groups and
their full unipotent subgroups.
We give our analysis in terms of the number of basic operations in the algebra 
$\E$:
addition, negation, multiplication, and testing equality.
Once again let $G$ be an $\F$-split reductive algebraic group, with
split maximal torus $T$, and Borel subgroup $B$ containing $T$.
Let $U$ be the unipotent radical of $B$.
Let $W=N_G(T)/T$ be the Weyl group and let $\Phi$ be the root system.
The reflection in $W$ corresponding to the root $\al$ is denoted $s_\al$.
Let $\Phi^+$ be the positive roots with respect to $B$.

First we give an analysis for element operations in $U(\E)$.
\begin{Thm}\label{T-analysisU}
Let $\F$ be a field and let $\E$ be a commutative unital $\F$-algebra.
Let $U$ be the full unipotent subgroup of a split reductive linear algebraic
group $G$ over $\F$.
Let $\ell$ be the semisimple rank of $G$.
Then there is a normal form for elements of $U(\E)$.
The word problem for elements in normal form requires $O(\ell^2)$ algebra operations, 
and multiplying or inverting them requires $O(\ell^3)$ algebra operations. 
\end{Thm}
\begin{proof}
The normal form is a collected word, so the timing for the word problem
follows from the fact that $N=|\Phi^+|$ is $O(\ell^2)$.
We can assume that $G$ is simple, since $U$ is a direct sum of the full unipotent
subgroups of the simple components of $G$.

If $G$ is classical, the formulas of Section~\ref{S-class} require $O(\ell^3)$
field operations.

If $G$ is exceptional, then $\ell$ is bounded.
In this case we use symbolic collection.  
The Hall polynomials of the full unipotent subgroup of a split reductive group
are independent of the algebra $\E$, since split reductive groups can be
constructed as $\Z$-schemes. 
So the number of algebra operations required for inversion or multiplication is
independent of the choice of $\E$. 
\end{proof}

In the rest of this section, we take $\E$ to be an extension field of $\F$
and we add inversion to the list of basic operations in $\E$.
We are primarily interested in computing in $U(\E)$ because it allows us to
compute in $G(\E)$ with the algorithms of \cite[Section~5]{CohenMurrayTaylor04}.
Recall that $G(\E)$ has a Steinberg presentation with generators $x_\al(a)$, for $\al\in\Phi$
and $a\in\E$; $n_\al$, for $\al\in\Phi$; and $t\in T(\E)$.
Note that the generator $x_r(a)$ of Section~\ref{S-collout} can be identified
with the generator $x_{\al_r}(a)$ of the Steinberg presentation.
Every element $g\in G(\E)$ can be written uniquely in Bruhat form:
$$
  g = ut\wdot u',
$$
for
\begin{itemize}
\item $u\in U(\E)$ stored as a collected word;
\item $t\in T(\E)$ stored as in \cite{CohenMurrayTaylor04};
\item $\wdot=n_{\al_1}\cdots n_{\al_m}$,  where $s_{\al_1}\cdots s_{\al_m}$ is a 
reduced expression for $w\in W$; and
\item $u'\in U_w(\E)$  as a collected word, where $U_w$ is the subgroup of $U$ 
generated by the terms 
$x_\al(a)$, for $\al$ in
$
  \Phi_w:=\{\al\in\Phi^+\mid \al w^{-1}\notin\Phi^+\}.
$
\end{itemize}
Given two elements in Bruhat form, we need to find the Bruhat form of their
product.
The usual element operations in $U(\E)$ are not sufficient for this
purpose.
There are two difficult steps, each of which requires a new operation in $U$.
We now describe these operations and show how to carry them out 
with the methods of the previous sections.

\subsection{Single-term separation}
One difficult step is multiplying $g = ut\wdot u'$ by $n_\al$ for some
$\al\in\Phi$.
This is achieved with Algorithm~3 of \cite{CohenMurrayTaylor04}, which uses the
following operation:
write $u'=\prod_{\be\in\Phi_w}x_\be(a_\be)$ in the form $x_\al(a_\al)v$ where
$v = \prod_{\be\in\Phi_w\setminus\{\al\}}x_\be(b_\be)$.
We call this operation \emph{single-term separation}.

This is easily done by collection: simply collect the term $x_\al(a_\al)$ to the
front of the product as in collection to the left, then put $v$ in the required
form with collection from the outside.
No extra terms of the form $x_\al(b)$ can appear in $v$ because only terms
corresponding to roots higher than $\al$ are created.
We can also do single term separation symbolically as in
Section~\ref{S-coll}.

Alternatively, for classical groups, we can compute $v$ as the product 
$x_\al(-a_\al)u'$ using the formulas of Section~\ref{S-class}.
If $\al=\al_{ij}$, then the only possible nonzero constants in $\varphi(\bfa)$ 
are $a_{ij}$, $a'_{ij}$, and $a''_i$.
Hence at most $O(\ell)$ of the formulas for $c_{ij}$ are nontrivial.
Each such formula has at most a constant number of nonzero terms. 
We now have:
\begin{Prop}\label{P-single}
Single-term separation in $U(\E)$ requires $O(\ell)$ field operations.
\end{Prop}
\begin{proof}
Use formulas for classical components and symbolic collection for exceptional
components.
\end{proof}

Note that, when both of the elements being multiplied are in Bruhat form,
Algorithm~3 of \cite{CohenMurrayTaylor04} only uses single-term separation for
$\al$ simple.
We have considered nonsimple roots as well, because they will be useful in the next
subsection. 

\subsection{Weyl separation}
The other difficult step for multiplication in $G$ is computing the product of 
$g\in G$ and $v\in U$.
Write $g$ in Bruhat form $ut\wdot u'$.
Then multiply $u'$ and $v$, and decompose the product into the
form  $v''v'$ where
$$ 
  v'' =\prod_{\al\in\Phi^+\setminus\Phi_w}x_\al(b_\al)\quad\text{and}
  \quad v'=\prod_{\al\in\Phi_w}x_\al(b_\al).
$$
We call this operation\/ \emph{Weyl separation}.
We now get the  Bruhat form
$$
  gv = [u(v'')^{\wdot^{-1}t^{-1}}]t\wdot v'
$$
where  $(v'')^{\wdot^{-1}t^{-1}}$ is in $U$ since $\al\in\Phi^+\setminus\Phi_w$
implies $\al w^{-1}$ is positive.

If we take the elements of $\Phi^+\setminus\Phi_w$ in an order compatible with 
height, followed by the elements of $\Phi_w$ in an order compatible with height,
we get a left-additive ordering on $\Phi^+$.
So the algorithms of Section~\ref{S-coll} can also be used for separation.
But note that $(v'')^{\wdot^{-1}t^{-1}}$ will need to be collected again, since
the image of the left additive ordering on $\Phi^+\setminus\Phi_w$ under
$w^{-1}$ need not be left additive.

We can also use collection from the outside for Weyl separation.
We need the following classification of additive orderings from \cite{Papi94}:
\begin{Thm}
Let $w$ be an element of the Weyl group $W$.
Let $s_{\be_1}\cdots s_{\be_m}$ be a reduced expression for $w$.
Then
$$
  %\be_N,\; \be_{N-1}s_{\be_N},\; \be_{N-2}s_{\be_{N-1}}s_{\be_N},\;  \dots,\; \be_1s_{\be_2}\cdots s_{\be_N}.
  \be_1s_{\be_2}\cdots s_{\be_N}, \dots,\; 
  \be_{N-2}s_{\be_{N-1}}s_{\be_N},\; \be_{N-1}s_{\be_N},\; \be_N
$$
is an additive ordering on $\Phi_w$.
All additive orderings on $\Phi_w$ arise from reduced expressions in this manner.
\end{Thm}
\noindent
Now let $w_0$ be the longest word in $W$ and fix a reduced expression
$s_{\al_1}\cdots s_{\al_N}$ for~$w_0$
(in practice, we use the lexicographically least reduced expression, but this is not
necessary).
We use the additive ordering on $\Phi^+$ corresponding to this reduced
expression as the fixed order for collection.
Now let $w$ be a Weyl group element.
If we restrict the fixed ordering to $\Phi_w$ we get an additive ordering, with
corresponding reduced expression $s_{\be_1}\dots s_{\be_m}=w$.
Similarly we restrict to get an ordering on 
$\Phi_{w_0w^{-1}} = (\Phi^+\setminus\Phi_w)w^{-1}$ 
and a corresponding 
reduced expression $s_{\ga_1}\dots s_{\ga_{N-m}}=w_0w^{-1}$.
Now
$w_0=s_{\ga_1}\dots s_{\ga_{N-m}}s_{\be_1}\dots s_{\be_m}$ is reduced.
The corresponding ordering is: our fixed ordering restricted to
$\Phi^+\setminus\Phi_w$ and transformed by $w$, followed by our fixed
ordering restricted to $\Phi_w$.
This is precisely the ordering we need for separation.

Finally we analyse Weyl separation:
\begin{Prop}\label{P-weyl}
Weyl separation in $U$ requires $O(\ell^3)$ field operations.
\end{Prop}
\begin{proof}
For classical components, we apply single-term separation for each root in
$\Phi_w$.
By Proposition~\ref{P-single}, this takes $O(N\ell)=O(\ell^3)$ operations.
For exceptional components  use  symbolic collection.
\end{proof}
\noindent In the exceptional case, this proposition assumes we have a system of 
symbolic-collection
polynomials for every Weyl element.
Although this is polynomial time, the memory required to store all these 
polynomials is prohibitive.
In practice, it is much faster to use
collection from the outside for Weyl separation in exceptional groups.

\subsection{Operations in reductive groups}
We now prove the following result on computation in $G$:
\begin{Thm}\label{T-analysis}
Let $\F$ be a field and let $\E$ be an extension of $\F$.
Let $G$ be a split reductive linear algebraic group over the field $\F$.
Let $\ell$ be the semisimple rank of $G$ and let $n$ be the reductive rank.
Then there is a normal form for elements of $G(\E)$.
The word problem for elements in normal form requires $O(n+\ell^2)$ field operations, 
and multiplying or inverting them requires $O(n\ell^2)$ field operations. 
\end{Thm}
\begin{proof}
We use the Bruhat decomposition to store $g\in G$  in the normal form
$g=uh\wdot u'$.
Here $u$, $u'$, and $\wdot$ are words of length at most $N$, 
while $h$ has length $n$.
Once again the timing for the word problem is clear.
Now $G$ is a central product of simple algebraic groups and a central
torus of dimension at most $n$.
Multiplying a toral element by an element in a simple component is
done as in Subsection~5.5 of~\cite{CohenMurrayTaylor04}, and takes time $O(n\ell^2)$.
So it suffices to show that multiplication and inversion in a simple group 
$G$ requires $O(\ell^3)$ operations.
The algorithms for multiplication and inversion given in
\cite{CohenMurrayTaylor04} require a constant number of multiplications or
inversions in $U(\E)$, together with a constant number of Weyl separations and at
most $O(\ell^2)$ single-term separations.
The theorem now follows from Theorem~\ref{T-analysisU},
Proposition~\ref{P-single}, and Proposition~\ref{P-weyl}.
\end{proof}
\noindent Theorem~\ref{T-main} is an immediate consequence of this result and the
fact that $\ell\le n$.

\section{Implementation and timings}\label{S-time}
A number of heuristic improvements are built into our implementations of the
algorithms described.
Most of them are either obvious or were suggested by our profiling of the code.
We restrict ourselves here to a brief description of the basic data
types used.
Representations of elements of the field $\F$ or algebra $\E$ are 
taken care of by the Magma computer algebra system \cite{BosmaCannonPlayoust97}.
Most of our code is written in traditional C \cite{KernighanRitchie88Book} and 
incorporated into the Magma core.
Less time-critical code is written in the Magma language itself.

A collected product $\prod_{r=1}^N x_r(a_r)$ is stored as a sequence
$[a_1,\dots,a_N]$.
While doing the collection,
we represent a term $x_r(a)$ as a pair 
$( r,a )$ of an integer and an element of $\E$.
Note that pairs $(r,0)$ are trivial -- they are always eliminated as soon as
they occur.
A word $\prod_{i=1}^M x_{r_i}(a_i)^{\ep_i}$
is represented as a doubly linked chain.
That is, every root element in the chain contains a reference to
its predecessor and successor, which is a null-reference if the 
element is the first (resp.\ last) in the chain:
\newcommand{\bla}{\stackrel{\rightarrow}{\leftarrow}}
$$
    \leftarrow \circ \leftrightarrows \circ \leftrightarrows \dots 
    \leftrightarrows  \circ \rightarrow
$$
We use this data structure because inserting and deleting terms
in the word when applying relations \eqref{E-r}--\eqref{E-rs}
can be done in constant time. 
For sequences, insertion and deletion would be more expensive, since the tail of
the sequence has to be copied in memory. 
The chain is doubly linked, since we need both the predecessor and the 
successor of a term in the word for the {\sc CollectSubword} functions.

In our tables we use the following abbreviations for collection algorithms:
\begin{center}
\begin{tabular}{lp{10cm}}
CTL:&  Collection to the left, Section~\ref{S-coll}.\\
CFL:&  Collection from the left, Section~\ref{S-coll}.\\
CFO:&  Collection from the outside, Section~\ref{S-collout}.\\
SCFL:& Symbolic collection from the left, Section~\ref{S-coll}.\\
SCFO:& Symbolic collection from the outside, Section~\ref{S-collout}.\\
\end{tabular}
\end{center}
We have two different implementations of the method of Section~\ref{S-class}:
\begin{description}
\item[\rm D]
Modified matrix multiplication. For given $\bfa,\bfb$, we use formulas of 
Section~\ref{S-class} to compute $\varphi(\bfa)$ and the
significant part of $\varphi(\bfb)$ 
(we do not use $b'_{ij}$ and, except in types
$\CB$ (even characteristic) and $\CC$, we do not use $b''_i$). 
Then the product of the two
matrices is computed using algorithms implemented in the Magma computer algebra
system \cite{BosmaCannonPlayoust97}. The resulting matrix agrees with
$\varphi(\bfc)$ in the entries $c_{ij}$ and in the entries $c''_i$ (where they
are needed). 
Thus we can recover  the product $\bfc=\bfa\bfb$ from the
resulting matrix by formulas of Section~\ref{S-class}.
\item[\rm SD]
Compute the polynomials of~(\ref{E-symb}) in Section~\ref{S-coll}, using the 
formulas instead of collection.
\end{description}
Method D outperforms SD in most cases, since 
asymptotically fast algorithms are used for matrix multiplication. But 
SD is faster for fields with very rapid blow-up of terms, 
such as multivariate rational function fields.
All timings were run on an AMD Opteron 150 Processor with 2393 MHz.

Table~\ref{T-preproc-times} gives times and memory consumption for
creating the reductive groups and precomputing all constants.
For symbolic algorithms, this also includes time taken to compute the
polynomials.
Note that all constants and polynomials are independent of the
field, and are computed on a per-root-datum basis. 
This means
that preprocessing time is nearly zero if a group with 
the same root datum has already been created in the same Magma session.
We used a workspace of $4$~gigabytes -- when this is insufficient we do not give a time
and write $>4GB$ in the memory column.
In columns D and SD, the constants are computed as they are needed and 
not stored in memory. 
\begin{sidewaystable}
%\begin{table}
{\scriptsize%\hspace*{-2cm}
\begin{tabular}{|l|rrrrrrr|rrrrrrr|}\hline
     		&\multicolumn{6}{l}{Time}                   && \multicolumn{6}{l}{Memory (MB)} & \\
 		& CTL & CFL & CFO & SCFL & SCFO & D  & SD   & CTL & CFL & CFO & SCFL & SCFO & D  & SD \\
\hline
$\CA_{ 10}(17)$ &   0.790 & 0.790 & 0.780 & 0.780 & 0.820 & 0.180 & 0.200    & 4.469 & 4.469 & 4.469 & 4.674 & 4.659 & 3.726 & 3.800   \\
$\CA_{ 20}(17)$ &   8.950 & 8.870 & 8.990 & 9.000 & 9.310 & 0.180 & 0.190    & 12.853 & 12.853 & 12.853 & 15.488 & 13.821 & 3.708 & 4.227   \\
$\CA_{ 30}(17)$ &   44.110 & 44.330 & 44.150 & 46.670 & 45.580 & 0.200 & 0.240    & 53.915 & 53.915 & 53.915 & 79.217 & 58.180 & 3.781 & 5.416   \\
$\CA_{100}(17)$ &   -- & -- & -- & -- & -- & 0.640 & 2.250    & $> 4GB$ & $> 4GB$ & $> 4GB$ & $> 4GB$ & $> 4GB$ & 8.139 & 52.477   \\
\hline
$\CB_{ 10}(17)$ &   4.490 & 4.320 & 4.480 & 4.390 & 4.580 & 0.180 & 0.180    & 4.731 & 4.731 & 4.731 & 6.219 & 6.377 & 3.693 & 3.992   \\
$\CB_{ 20}(17)$ &   68.700 & 68.870 & 68.570 & 72.270 & 69.750 & 0.190 & 0.250    & 39.072 & 39.072 & 39.072 & 83.714 & 45.900 & 3.710 & 7.951   \\
$\CB_{ 30}(17)$ &   357.170 & 357.060 & 355.760 & 437.320 & 365.790 & 0.200 & 0.430    & 177.452 & 177.452 & 177.452 & 613.679 & 426.778 & 3.790 & 14.250   \\
$\CB_{100}(17)$ &   -- & -- & -- & -- & -- & 0.620 & 14.380    & $> 4GB$ & $> 4GB$ & $> 4GB$ & $> 4GB$ & $> 4GB$ & 8.160 & 445.203   \\
\hline
$\CC_{ 10}(17)$ &   4.460 & 4.400 & 4.560 & 4.380 & 4.530 & 0.180 & 0.200    & 4.730 & 4.730 & 4.730 & 6.218 & 5.549 & 3.693 & 3.992   \\
$\CC_{ 20}(17)$ &   68.730 & 68.680 & 68.800 & 72.720 & 69.870 & 0.190 & 0.250    & 39.078 & 39.078 & 39.078 & 83.721 & 65.053 & 3.710 & 7.951   \\
$\CC_{ 30}(17)$ &   354.660 & 357.230 & 354.970 & 436.860 & 363.670 & 0.190 & 0.440    & 177.449 & 177.449 & 177.449 & 613.676 & 285.094 & 3.790 & 14.250   \\
$\CC_{100}(17)$ &   -- & -- & -- & -- & -- & 0.620 & 16.490    & $> 4GB$ & $> 4GB$ & $> 4GB$ & $> 4GB$ & $> 4GB$ & 8.160 & 445.205   \\
\hline
$\CD_{ 10}(17)$ &   1.820 & 1.820 & 1.710 & 1.820 & 1.870 & 0.190 & 0.160    & 4.344 & 4.344 & 4.344 & 5.198 & 4.895 & 3.692 & 3.986   \\
$\CD_{ 20}(17)$ &   28.540 & 28.400 & 28.400 & 31.640 & 29.910 & 0.190 & 0.270    & 34.589 & 34.589 & 34.589 & 79.227 & 54.305 & 3.708 & 7.911   \\
$\CD_{ 30}(17)$ &   153.440 & 153.080 & 152.900 & 220.920 & 159.850 & 0.200 & 0.400    & 167.857 & 167.857 & 167.857 & 604.077 & 249.470 & 3.785 & 14.242   \\
$\CD_{100}(17)$ &   -- & -- & -- & -- & -- & 0.640 & 14.180    & $> 4GB$ & $> 4GB$ & $> 4GB$ & $> 4GB$ & $> 4GB$ & 8.158 & 445.176   \\
\hline
$\CG_{  2}(17)$ &   0.210 & 0.190 & 0.190 & 0.190 & 0.200 &  --   &  --      & 3.580 & 3.580 & 3.580 & 3.580 & 3.580 &  --   &  --     \\
$\CF_{  4}(17)$ &   0.410 & 0.440 & 0.400 & 0.400 & 0.430 &  --   &  --      & 3.764 & 3.764 & 3.764 & 3.859 & 3.764 &  --   &  --     \\
$\CE_{  6}(17)$ &   0.440 & 0.450 & 0.430 & 0.440 & 0.460 &  --   &  --      & 4.532 & 4.532 & 4.532 & 4.532 & 4.532 &  --   &  --     \\
$\CE_{  7}(17)$ &   0.990 & 1.000 & 0.980 & 1.010 & 0.980 &  --   &  --      & 4.726 & 4.726 & 4.726 & 4.305 & 3.870 &  --   &  --     \\
$\CE_{  8}(17)$ &   3.130 & 3.110 & 3.100 & 3.180 & 3.230 &  --   &  --      & 6.088 & 6.088 & 6.088 & 14.175 & 7.376 &  --   &  --     \\
\hline

\end{tabular}}\\
\caption{Time and memory consumption for preprocessing}
\label{T-preproc-times}
%\end{table}
\end{sidewaystable}

Table~\ref{T-unip-times} gives average times for multiplying and inverting 
random elements of
full unipotent groups over the field with 17 elements.
The average is taken over $100$ multiplications. 
The same random elements are used for different algorithms.
If a single multiplication required more than $2$~gigabytes of memory, we write 
$>2GB$
instead of a time.
We did not attempt those cases where the preprocessing took more than
$4$~gigabytes of
memory.
\begin{sidewaystable}
%\begin{table}
{\scriptsize%\hspace*{-2cm}
\begin{tabular}{|l|rrrrrrr|rrrrrrr|}\hline
Group		& \multicolumn{6}{l}{Multiply}              &&\multicolumn{6}{l}{Invert} & \\
 		& CTL & CFL & CFO & SCFL & SCFO & D  & SD   & CTL & CFL & CFO & SCFL & SCFO & D  & SD \\
\hline
$\CA_{ 10}(17)$ &   0.009 & 0.006 & 0.006 & 0.007 & 0.006 & 0.006 & 0.006    & 0.003 & 0.002 & 0.002 & 0.003 & 0.002 & 0.001 & 0.002   \\
$\CA_{ 20}(17)$ &   1.511 & 0.058 & 0.030 & 0.174 & 0.031 & 0.020 & 0.032    & 1.051 & 0.042 & 0.012 & 0.157 & 0.019 & 0.003 & 0.017   \\
$\CA_{ 30}(17)$ &   114.543 & 0.768 & 0.111 & 3.069 & 0.334 & 0.046 & 0.182    & 84.299 & 0.734 & 0.063 & 3.052 & 0.326 & 0.007 & 0.148   \\
$\CA_{100}(17)$ &   -- & -- & -- & -- & -- & 0.854 & 51.744    & -- & -- & -- & -- & -- & 0.071 & 52.425   \\
\hline
$\CB_{ 10}(17)$ &   0.101 & 0.016 & 0.013 & 0.045 & 0.038 & 0.012 & 0.014    & 0.069 & 0.009 & 0.005 & 0.039 & 0.032 & 0.004 & 0.007   \\
$\CB_{ 20}(17)$ &   280.288 & 1.144 & 0.180 & 5.134 & 0.614 & 0.050 & 0.321    & 209.010 & 1.105 & 0.136 & 5.474 & 0.622 & 0.017 & 0.305   \\
$\CB_{ 30}(17)$ &   $> 2GB$ & 25.472 & 1.412 & 107.734 & 49.190 & 0.122 & 2.015    & $> 2GB$ & 25.419 & 1.240 & 115.646 & 49.338 & 0.045 & 2.024   \\
$\CB_{100}(17)$ &   -- & -- & -- & -- & -- & 2.728 & 1025.957    & -- & -- & -- & -- & -- & 1.115 & 1048.822   \\
\hline
$\CC_{ 10}(17)$ &   0.093 & 0.016 & 0.013 & 0.044 & 0.016 & 0.013 & 0.014    & 0.065 & 0.009 & 0.005 & 0.038 & 0.009 & 0.004 & 0.007   \\
$\CC_{ 20}(17)$ &   266.161 & 1.133 & 0.178 & 4.999 & 2.306 & 0.050 & 0.319    & 194.365 & 1.106 & 0.119 & 5.411 & 2.316 & 0.017 & 0.300   \\
$\CC_{ 30}(17)$ &   $> 2GB$ & 27.562 & 1.443 & 113.089 & 21.446 & 0.125 & 2.099    & $> 2GB$ & 27.656 & 1.097 & 123.097 & 21.609 & 0.047 & 2.098   \\
$\CC_{100}(17)$ &   -- & -- & -- & -- & -- & 2.760 & 1046.944    & -- & -- & -- & -- & -- & 1.161 & 1054.659   \\
\hline
$\CD_{ 10}(17)$ &   0.062 & 0.014 & 0.011 & 0.020 & 0.012 & 0.011 & 0.012    & 0.040 & 0.007 & 0.004 & 0.014 & 0.006 & 0.004 & 0.006   \\
$\CD_{ 20}(17)$ &   195.012 & 0.887 & 0.144 & 4.470 & 1.621 & 0.046 & 0.293    & 141.648 & 0.874 & 0.096 & 4.713 & 1.619 & 0.016 & 0.275   \\
$\CD_{ 30}(17)$ &   $> 2GB$ & 22.057 & 1.181 & 102.313 & 15.470 & 0.121 & 1.963    & $> 2GB$ & 22.229 & 0.931 & 108.976 & 15.613 & 0.045 & 1.969   \\
$\CD_{100}(17)$ &   -- & -- & -- & -- & -- & 2.631 & 1038.942    & -- & -- & -- & -- & -- & 1.093 & 1044.249   \\
\hline
$\CG_{  2}(17)$ &   0.001 & 0.001 & 0.001 & 0.001 & 0.001 &  --   &  --      & 0.000 & 0.000 & 0.000 & 0.000 & 0.000 &  --   &  --     \\
$\CF_{  4}(17)$ &   0.004 & 0.003 & 0.003 & 0.003 & 0.003 &  --   &  --      & 0.001 & 0.001 & 0.001 & 0.001 & 0.001 &  --   &  --     \\
$\CE_{  6}(17)$ &   0.006 & 0.004 & 0.004 & 0.005 & 0.004 &  --   &  --      & 0.002 & 0.001 & 0.001 & 0.002 & 0.001 &  --   &  --     \\
$\CE_{  7}(17)$ &   0.042 & 0.009 & 0.007 & 0.014 & 0.008 &  --   &  --      & 0.029 & 0.004 & 0.002 & 0.009 & 0.004 &  --   &  --     \\
$\CE_{  8}(17)$ &   4.924 & 0.053 & 0.016 & 0.309 & 0.032 &  --   &  --      & 3.966 & 0.044 & 0.008 & 0.292 & 0.027 &  --   &  --     \\
\hline

\end{tabular}}\\
\caption{Average time to multiply random elements of the full unipotent group}
\label{T-unip-times}
%\end{table}
\end{sidewaystable}

Table~\ref{T-red-times} gives similar times for multiplying and inverting random elements of
the reductive group itself.  
Each such operation involves a number of collections.
Computing random elements in a reductive group can be time consuming, but this
is not included in our timings.
\begin{sidewaystable}
%\begin{table}
{\scriptsize%\hspace*{-2cm}
\begin{tabular}{|l|rrrrrrr|rrrrrrr|}\hline
Group		& \multicolumn{6}{l}{Multiply}              &&\multicolumn{6}{l}{Invert} & \\
 		& CTL & CFL & CFO & SCFL & SCFO & D  & SD   & CTL & CFL & CFO & SCFL & SCFO & D  & SD \\
\hline
$\CA_{ 10}(17)$ &   0.123 & 0.119 & 0.118 & 0.121 & 0.117 & 0.288 & 0.295    & 0.163 & 0.160 & 0.162 & 0.163 & 0.160 & 0.222 & 0.220   \\
$\CA_{ 20}(17)$ &   1.359 & 0.729 & 0.662 & 1.293 & 0.695 & 2.063 & 2.129    & 2.268 & 1.147 & 1.080 & 1.546 & 1.112 & 1.982 & 2.031   \\
$\CA_{ 30}(17)$ &   120.524 & 10.798 & 7.887 & 21.833 & 8.936 & 29.141 & 29.752    & 104.102 & 13.639 & 12.655 & 21.838 & 13.484 & 21.959 & 22.627   \\
$\CA_{100}(17)$ &   -- & -- & -- & -- & -- &       &          & -- & -- & -- & -- & -- &       &         \\
\hline
$\CB_{ 10}(17)$ &   0.423 & 0.325 & 0.304 & 0.453 & 0.414 & 1.353 & 1.359    & 0.561 & 0.495 & 0.486 & 0.592 & 0.570 & 0.860 & 0.867   \\
$\CB_{ 20}(17)$ &   342.015 & 9.964 & 5.184 & 28.542 & 7.234 & 27.707 & 29.383    & 237.331 & 10.345 & 8.721 & 24.834 & 10.331 & 18.471 & 19.728   \\
$\CB_{ 30}(17)$ &   $> 2GB$ & 170.359 & 34.042 & 560.664 & 228.163 & 190.445 & 198.912    & $> 2GB$ & 92.252 & 51.170 & 397.654 & 197.548 & 140.195 & 146.980   \\
$\CB_{100}(17)$ &   -- & -- & -- & -- & -- &       &          & -- & -- & -- & -- & -- &       &         \\
\hline
$\CC_{ 10}(17)$ &   0.430 & 0.329 & 0.314 & 0.455 & 0.334 & 1.420 & 1.426    & 0.556 & 0.488 & 0.483 & 0.588 & 0.502 & 0.893 & 0.905   \\
$\CC_{ 20}(17)$ &   347.298 & 10.341 & 5.380 & 28.422 & 14.166 & 27.950 & 29.116    & 225.542 & 10.033 & 8.517 & 24.409 & 15.189 & 18.399 & 19.387   \\
$\CC_{ 30}(17)$ &   $> 2GB$ & 174.191 & 33.764 & 580.168 & 116.204 & 181.660 & 190.782    & $> 2GB$ & 99.784 & 54.166 & 414.583 & 115.765 & 138.295 & 145.607   \\
$\CC_{100}(17)$ &   -- & -- & -- & -- & -- &       &          & -- & -- & -- & -- & -- &       &         \\
\hline
$\CD_{ 10}(17)$ &   0.344 & 0.280 & 0.268 & 0.317 & 0.278 & 1.170 & 1.184    & 0.454 & 0.412 & 0.409 & 0.441 & 0.416 & 0.727 & 0.734   \\
$\CD_{ 20}(17)$ &   225.229 & 8.183 & 4.606 & 24.562 & 10.825 & 23.712 & 24.692    & 166.790 & 8.538 & 7.340 & 21.308 & 12.015 & 15.674 & 16.571   \\
$\CD_{ 30}(17)$ &   $> 2GB$ & 147.254 & 31.131 & 521.051 & 90.824 & 174.995 & 185.490    & $> 2GB$ & 82.122 & 47.233 & 373.211 & 94.738 & 135.237 & 142.438   \\
$\CD_{100}(17)$ &   -- & -- & -- & -- & -- &       &          & -- & -- & -- & -- & -- &       &         \\
\hline
$\CG_{  2}(17)$ &   0.008 & 0.008 & 0.006 & 0.008 & 0.008 &  --   &  --      & 0.006 & 0.006 & 0.006 & 0.006 & 0.006 &  --   &  --     \\
$\CF_{  4}(17)$ &   0.080 & 0.076 & 0.075 & 0.077 & 0.074 &  --   &  --      & 0.060 & 0.060 & 0.059 & 0.059 & 0.060 &  --   &  --     \\
$\CE_{  6}(17)$ &   0.172 & 0.166 & 0.158 & 0.165 & 0.158 &  --   &  --      & 0.128 & 0.129 & 0.124 & 0.127 & 0.126 &  --   &  --     \\
$\CE_{  7}(17)$ &   0.774 & 0.508 & 0.451 & 0.518 & 0.455 &  --   &  --      & 0.419 & 0.369 & 0.363 & 0.381 & 0.366 &  --   &  --     \\
$\CE_{  8}(17)$ &   52.481 & 3.098 & 1.566 & 3.815 & 1.624 &  --   &  --      & 8.781 & 1.396 & 1.288 & 1.974 & 1.353 &  --   &  --     \\
\hline

\end{tabular}}\\
\caption{Average time to multiply random elements of the reductive group}
\label{T-red-times}
%\end{table}
\end{sidewaystable}

Table~\ref{T-fld-times} gives average times for multiplying and inverting random elements of
full unipotent groups of reductive groups over different fields.
Over the field of rational numbers, the random field elements are chosen by
taking a random 
numerator and a random denominator of size up to the given number of bits and a random sign.
Similar random elements were used for the Gaussian integers $\Q(i)$ and for
$\Q(p)$, which is the splitting field of a random irreducible polynomial of
degree 6 with integral coefficients in the range 1 to 10.
The field $R$ is the multivariate rational function field over $\Q$ with 10
variables. Random field elements over $R$ were taken to be random invariates.
In $\CB_{20}(R)$ the coefficient blowup is so large that over $2$~gigabytes of 
memory was needed in some cases (see entries in the table).
%
%\begin{sidewaystable}
\begin{table}
{\scriptsize
\begin{tabular}{|l|rrrr|rrrr|}\hline
Group		& \multicolumn{3}{l}{Multiply}   &&\multicolumn{3}{l}{Invert} & \\
 		& CFO & SCFO & D  & SD           & CFO & SCFO & D  & SD \\
\hline
$\CA_{ 20}( 2)$ &   0.021 & 0.030 & 0.020 & 0.031     & 0.005 & 0.017 & 0.003 & 0.015   \\
$\CA_{ 20}(17)$ &   0.029 & 0.031 & 0.020 & 0.032     & 0.012 & 0.018 & 0.003 & 0.017   \\
$\CA_{ 20}(\Q)$, \hfill$ 32$ bits &   0.045 & 0.038 & 0.024 & 0.038     & 0.051 & 0.109 & 0.213 & 0.135   \\
$\CA_{ 20}(\Q)$, \hfill$ 64$ bits &   0.049 & 0.041 & 0.028 & 0.041     & 0.086 & 0.195 & 0.736 & 0.265   \\
$\CA_{ 20}(\Q)$, \hfill$128$ bits &   0.056 & 0.047 & 0.039 & 0.047     & 0.177 & 0.394 & 2.877 & 0.611   \\
$\CA_{ 20}(\Q(i))$, \hfill$ 32$ bits &   0.074 & 0.052 & 0.047 & 0.049     & 0.080 & 0.143 & 0.048 & 0.117   \\
$\CA_{ 20}(\Q(p))$, \hfill$ 32$ bits &   0.108 & 0.069 & 0.062 & 0.062     & 0.122 & 0.265 & 0.068 & 0.214   \\
$\CA_{ 20}( R)$ &   0.049 & 0.038 & 0.025 & 0.038     & 6.868 & 2.562 & 0.594 & 1.609   \\
\hline
$\CB_{ 20}( 2)$ &   0.056 & 0.589 & 0.047 & 0.303     & 0.023 & 0.590 & 0.015 & 0.282   \\
$\CB_{ 20}(17)$ &   0.174 & 0.592 & 0.049 & 0.309     & 0.125 & 0.596 & 0.017 & 0.289   \\
$\CB_{ 20}(\Q)$, \hfill$ 32$ bits &   0.420 & 0.811 & 0.254 & 1.901     & 0.986 & 2.630 & 2.625 & 5.169   \\
$\CB_{ 20}(\Q)$, \hfill$ 64$ bits &   0.534 & 0.960 & 0.630 & 3.594     & 2.205 & 5.111 & 11.937 & 12.567   \\
$\CB_{ 20}(\Q)$, \hfill$128$ bits &   0.798 & 1.269 & 1.865 & 8.058     & 5.575 & 11.475 & 50.818 & 34.318   \\
$\CB_{ 20}(\Q(i))$, \hfill$ 32$ bits &   0.680 & 0.972 & 1.078 & 2.124     & 1.270 & 2.654 & 1.344 & 3.874   \\
$\CB_{ 20}(\Q(p))$, \hfill$ 32$ bits &   1.092 & 1.362 & 1.410 & 3.864     & 2.276 & 4.542 & 2.053 & 7.212   \\
$\CB_{ 20}( R)$ &   0.589 & 0.777 & $>2GB$ & 31.884     & $>2GB$ & $>2GB$ & $>2GB$ & $>2GB$   \\
\hline
$\CE_{  8}( 2)$ &   0.012 & 0.028 &  --   &  --       & 0.003 & 0.022 &  --   &  --     \\
$\CE_{  8}(17)$ &   0.016 & 0.029 &  --   &  --       & 0.007 & 0.023 &  --   &  --     \\
$\CE_{  8}(\Q)$, \hfill$ 32$ bits &   0.047 & 0.125 &  --   &  --       & 0.061 & 0.283 &  --   &  --     \\
$\CE_{  8}(\Q)$, \hfill$ 64$ bits &   0.075 & 0.206 &  --   &  --       & 0.118 & 0.526 &  --   &  --     \\
$\CE_{  8}(\Q)$, \hfill$128$ bits &   0.141 & 0.390 &  --   &  --       & 0.262 & 1.099 &  --   &  --     \\
$\CE_{  8}(\Q(i))$, \hfill$ 32$ bits &   0.077 & 0.162 &  --   &  --       & 0.094 & 0.310 &  --   &  --     \\
$\CE_{  8}(\Q(p))$, \hfill$ 32$ bits &   0.113 & 0.291 &  --   &  --       & 0.146 & 0.621 &  --   &  --     \\
$\CE_{  8}( R)$ &   0.315 & 0.349 &  --   &  --       & 0.904 & 1.166 &  --   &  --     \\
\hline

\end{tabular}}
\vspace*{1mm}
\caption{Operations for random elements of the full unipotent group over different fields}
\label{T-fld-times}
\end{table}
%\end{sidewaystable}

Finally, in Table~\ref{T-Halldegs}, we compare the total degrees of the
polynomials used for different kinds of symbolic collection.
The last column contains 
$\lim_{\ell\to\infty} \mathrm{avg}\{ \deg(p) : p \in \mathcal{P} \}$,
where $\mathcal{P}$ is the set of polynomials used for symbolic collection.
This goes a long way towards explaining why collection from the outside works so
well. 
Since the polynomials are multivariate in $2N = 2|\Phi^+|$ variables,
we still can have large polynomials. To have a very rough idea of the size,
we printed the largest of the polynomials in type $\CB_{15}$ as a string
and measured its size in bytes. Using collection from the outside, the size is 
$9226$~bytes; using collection from the left, the size is about $297$~megabytes.
\begin{table}

\addtolength{\extrarowheight}{1pt}
{\scriptsize
\begin{tabular}{|l|c|rr|c|cr|}
\hline
&\multicolumn{2}{c}{CFL}&&\multicolumn{2}{c}{CFO}& \\
\hline
& \quad max \quad\ & \quad avg\qquad\ \ && \quad max \quad\  & \  lim avg  &\\
\hline
$\CA_\ell$ &  $\ell$   &  $(\ell+2)/3$                        &&  $2$  &  $2$ &\\
$\CB_\ell$ & $2\ell-1$ &  $(2\ell+\frac32-\frac1{2\ell})/3$   &&  $4$  &  $4$ &\\
$\CC_\ell$ & $2\ell-1$ &  $(2\ell+\frac32-\frac1{2\ell})/3$   &&  $3$  &  $3$ &\\
$\CD_\ell$ & $2\ell-3$ &  $(2\ell-1)/3$                       &&  $3$  &  $3$ &\\
\hline
\end{tabular}}
\vspace*{1mm}
\caption{Total degrees of Hall polynomials}\label{T-Halldegs}
\end{table}

\bibliographystyle{amsalpha}
\bibliography{references}

\end{document}